\documentclass[graybox]{svmult}

\usepackage{makeidx}         % allows index generation
\usepackage{graphicx}        % standard LaTeX graphics tool
\usepackage{multicol}        % used for the two-column index
\usepackage[bottom]{footmisc}% places footnotes at page bottom

\usepackage{newtxtext}       %
\usepackage{newtxmath}       % selects Times Roman as basic font

\usepackage{hyperref}
\usepackage{listings}
\usepackage{mathtools}

\makeindex             % used for the subject index
\renewcommand{\emph}[1]{\textsl{#1}}

\begin{document}

\title*{Teaching Programming for Mathematical Scientists}
\author{Jack Betteridge and Eunice Y.~S.~Chan and Robert M.~Corless and~James~H.~Davenport and James Grant }
\institute{Jack Betteridge \at The University of Bath, Bath, England \email{jdb55@bath.ac.uk}
\and Eunice Y.~S.~Chan \at Western University, London, Canada \email{echan295@uwo.ca}
\and Robert M.~Corless \at Western University, London, Canada \email{rcorless@uwo.ca}
\and \hbox{James H.~Davenport} \at The University of Bath, Bath, England \email{masjhd@bath.ac.uk}
\and James Grant \at The University of Bath, Bath, England \email{rjg20@bath.ac.uk}
}
\authorrunning{Betteridge \emph{et al}}
\maketitle

%\abstract*{Each chapter should be preceded by an abstract (no more than 200 words) that summarizes the content. The abstract will appear \textit{online} at \url{www.SpringerLink.com} and be available with unrestricted access. This allows unregistered users to read the abstract as a teaser for the complete chapter.
%Please use the 'starred' version of the \texttt{abstract} command for typesetting the text of the online abstracts (cf. source file of this chapter template \texttt{abstract}) and include them with the source files of your manuscript. Use the plain \texttt{abstract} command if the abstract is also to appear in the printed version of the book.}

\abstract{Over the past thirty years or so the authors have been teaching various \emph{programming for mathematics} courses at our respective Universities, as well as incorporating computer algebra and numerical computation into traditional mathematics courses.  These activities are, in some important ways, natural precursors to the use of Artificial Intelligence in Mathematics Education.
This paper reflects on some of our course designs and experiences and is therefore a mix of theory and practice.
Underlying both is a clear recognition of the value of computer programming for mathematics education. We use this theory and practice to suggest good techniques for and to raise questions about the use of AI in Mathematics Education. }

\section{Background}
\label{sec:background}
This volume is part of the fast-growing literature in a relatively new field---being only about thirty years old---namely Artificial Intelligence for Education (AIEd).  The survey~\cite{luckin2016intelligence} gives in lay language a concise overview of the field and advocates for its ambitious goals. For a well-researched discussion of an opposing view and of the limitations of Artificial Intelligence (AI) see~\cite{broussard2018artificial}. This present paper is concerned with AI in mathematics education in two senses: first, symbolic computation tools were themselves among the earliest and most successful pieces of AI to arise out of the original MIT labs already in the sixties\footnote{For example \cite{Slagle1963}, which took a `Good Old-Fashioned Artificial Intelligence (GOFAI)'' approach, and concluded ``The solution of a symbolic integration problem by a commercially available computer is far cheaper and faster than by man''. Of course this was from the era when people still believed in GOFAI.  We are grappling with different problems today, using much more powerful tools.  Yet some important things can be learned by looking at the effects of the simpler and older tools. The riposte to \cite{Slagle1963} was the development of Computer Algebra \cite{Davenport2018a} as a separate discipline.}, and have had a significant impact on mathematics education.  This impact is still changing the field of mathematics education, especially as the tools evolve~\cite{kovacs2017geogebra}.  Second, and we believe more important, the existence of these tools, and similarly the existence of other AI tools, has profoundly changed the affordances of mathematics and therefore \emph{should change the content of mathematics courses, not just their presentation methods}~\cite{corless2004computer}.  That last paper introduced the phrase ``Computer-Mediated Thinking,'' by which was meant an amplification of human capability by use of computer.  In~\cite{hegedus2017uses} we find this idea beautifully articulated and set in the evolving sequence of human methods for mediating their thinking: symbolic marks on bone, through language, to symbols on paper, to where we are today.  One of our theses is that people need to be given opportunities to learn how best to use these tools.
%, and how well the tools are taught affects how well the people think.

This paper reflects our experiences in changing mathematical course syllabi to reflect these new affordances, and may serve as a reference point in discussing future curricular changes (and what should not change) owing to the ever-shifting technological ground on which we work.  Our methodology is to consider mathematics education and computer programming together.  Our thesis is that the effect of computational tools, including AI, is greater and more beneficial if students are taught how to use the tools effectively and even how to create their own.

The connection between mathematics and computer programming is widely recognized and profound.  Indeed, most members of the general public will (if they think about it) simply assume that all mathematicians can, and do, program computers.  When mathematics is used instrumentally in science, as opposed to purely for aesthetic mathematical goals, this is in fact nearly universally true.  This is because computational tools are ubiquitous in the mathematical sciences.  Such tools are nowadays becoming increasingly accepted in pure mathematics, as well: see e.g.~\cite{borwein2009computer}. Modern tools even approach that most central aspect of pure mathematics, the notion of mathematical proof.  See~\cite{richard2019issues} and its many references for a careful and nuanced discussion of the notion of proof in a modern technologically-assisted environment, and the implications for mathematical education.

One lesson for educators is that we \emph{must} teach students in the mathematical sciences how to use computational tools responsibly. We earlier said that the public assumes that all mathematicians can program; with a similar justification, many students assume that they themselves can, too. But programming computers \emph{well} (or even just \emph{using} them well) is a different story.  The interesting thing for this paper is that learning to use computers well is a very effective way to learn mathematics well: by teaching programming, we can teach people to be better mathematicians and mathematical scientists.

We used the word ``must,'' above: we \emph{must} teach students how to $\ldots$.  Why ``must''?  For what reason?  We contend that this is the \emph{ethical} thing to do, in order to prepare our students as best we can to be functioning and thinking adults.  This is more than just preparation for a job: we are aiming at \emph{eudaemonia} here~\cite{flanagan2009really}.  This observation has significant implications for the current revolution in AI-assisted teaching.  We will return to this observation after discussing our experiences.

Our experience includes teaching programming to mathematical scientists and engineers through several eras of ``new technology,'' as they have flowed and ebbed. Our early teaching experience includes the use of computer-algebra capable calculators to teach engineering mathematics\footnote{We had intended to give the reference~\cite{Rosati:1992:Evaluation} for this; however, that journal seems to have disappeared and we can find no trace of it on the Web, which is a kind of testimony to ephemerality.  Some of the lessons of that article were specific to the calculator, which was \emph{too advanced} for its era and would be disallowed in schools today.  We shall not much discuss the current discouragingly restricted state of the use of calculators in schools hereafter.}; calculators were a good solution at the time because we could not then count on every student having their own computer (and smartphones were yet a distant technological gleam).  Some of the lessons we learned then are still valid, however: in particular, we learned that we must teach students that they are \emph{responsible} for the results obtained by computation, that they \emph{ought to know} when the results were reliable and when not, and that they should \emph{understand the limits of computation} (chiefly, understand both complexity of computation and numerical stability of floating-point computation; today we might add that generic algebraic results are not always valid under specialization, as in \cite{CamargosCoutoetal2020a}).  These lessons are invariant under shifts in technology, and become particularly pertinent when AI enters the picture.

Speaking of shifts, see~\cite{kahan1983mathematics} (``Mathematics written in Sand'') for an early attack on teaching those lessons, in what was then a purely numerical environment. A relevant quotation from that work is
\begin{quotation}
Rather than have to copy the received word,
students are entitled to experiment with
mathematical phenomena, discover more of them,
and then read how our predecessors discovered
even more. Students need inexpensive apparatus
analogous to the instruments and glassware in
Physics and Chemistry laboratories, but designed
to combat the drudgery that inhibits exploration.
\hfill---William Kahan, p.~5 \emph{loc cit.}
\end{quotation}

Teaching these principles in a way that the student can absorb them is a significant curricular goal, and room must be made for this goal in the mathematical syllabus.  This means that some things that are in that already overfull syllabus must be jettisoned.  In~\cite{corless1997scientific} and again in~\cite{corless2004computer} some of us claim that \emph{convergence tests for infinite series} should be among the first to go.  Needless to say, this is a radical proposal and not likely to attain universal adoption without a significant shift in policy; nevertheless, if not this, then what else?  Clearly \emph{something} has to go, to make room!

Curricular shifts are the norm, over time.  For instance, spherical trigonometry is no longer taught as part of the standard engineering mathematics curriculum; nor are graphical techniques for solving algebraic equations (which formerly were part of the \emph{drafting} curriculum, itself taken over by CAD).  Special functions are now taught as a mere rump of what they were, once.  Euclidean geometry has been almost completely removed from the high-school curriculum. Many of these changes happen ``by accident'' or for other, non-pedagogical, reasons; moreover it seems clear that removing Euclidean geometry has had a deleterious effect on the teaching of logic and proof, which was likely unintended.

We have found (and will detail some of our evidence for this below) that teaching \emph{programming} remediates some of these ill effects.  By learning to program, the student will in effect learn how to prove. If nothing else, learning to program may motivate the student wanting to \emph{prove the program correct}. This leads into the modern disciplines of Software Engineering and of Software Validation; not to mention Uncertainty Quantification. Of course there are some truly difficult problems hiding in this innocent-seeming suggestion: but there are advantages and benefits even to such intractable problems.

We will begin by discussing the teaching of Numerical Analysis and programming, in what is almost a traditional curriculum.  We will see some seeds of curriculum change in response to computational tools already in this pre-AI subject.

\section{Introduction to Numerical Analysis}
\label{sec:numericalanalysis}
The related disciplines of ``Numerical Methods,'' ``Numerical Analysis,'' ``Scientific Computing,'' and ``Computational Science'' need little introduction or justification nowadays (they could perhaps benefit from disambiguation).  Many undergraduate science and engineering degrees will grudgingly leave room for one programming course if it is called by one of those names.  Since this is typically the first course where the student has to actually \emph{use} the mathematical disciplines of linear algebra and calculus (and use them \emph{together}) there really isn't much room in such a course to teach good programming.  Indeed many students are appalled to learn that the techniques of \emph{real analysis}, itself a feared course, make numerical analysis intelligible.  In this minimal environment (at Western the course occupies $13$ weeks, with three hours of lecture and one\footnote{Students were enrolled in one of three tutorial hours, but often went to all three hours.} hour of tutorial per week) we endeavoured to teach the following:

\begin{enumerate}
  \item The basics of numerical analysis: \emph{backward error} and \emph{conditioning}
  \item How to write simple computer programs: conditionals, loops, vectors, recursion
  \item The elements of programming style: readability, good naming conventions, the use of comments
  \item Several important numerical algorithms: matrix factoring, polynomial approximation; solving IVP for ODE
  \item How to work in teams and to \emph{communicate mathematics}.
\end{enumerate}

The students also had some things to \emph{unlearn}: about the value of exact answers, or about the value of some algorithms that they had been taught to perform by hand, such as Cramer's Rule for solving linear systems of equations, for instance.

Western has good students, with an entering average amongst the highest in the country. By and large the students did well on these tasks.  But they had to work, in order to do well.  The trick was to get them to do the work.

\subsection{Choice of Programming Language}
\label{subsec:Choice}
We used Matlab.  This choice was controversial: some of our colleagues wanted us to teach C or C++ because, ultimately for large Computational Science applications, the speed of these compiled languages is necessary.  However, for the goals listed above, we think that Matlab is quite suitable; moreover, Matlab is a useful scientific language in its own right because \emph{development time} is minimized by programming in a high-level language first~\cite{wilson2014best}, and because of that Matlab is very widely used.

Other colleagues wanted us to use an open-source language such as Python, and this is quite attractive and Python may indeed eventually displace Matlab in this teaching role, for several reasons.  But as of this moment in time, Matlab retains some advantages in installed methods for solving several important problems and in particular its sparse matrix methods are very hard to beat.

We also used the computer algebra language Maple, on occasion: for comparison with exact numerical results, and for program generation.  Matlab's Symbolic Toolbox is quite capable but we preferred to separate symbolic computation from numeric computation for the purposes of the course.

\subsection{Pedagogical Methods}
We used \emph{Active Learning}, of course.  By now, the evidence in its favour is so strong as to indicate that \emph{not} using active learning is academically irresponsible~\cite{handelsman2004scientific,freeman2014active}.  However, using active learning techniques in an 8:30am lecture hall for 90 or so students in a course that is overfull of material is a challenge. To take only the simplest techniques up here, we first talk about \emph{Reading Memos}~\cite{smith1995teaching}.

We gave credit---five percent of the student's final mark---for simply \emph{handing in a Reading Memo}, that is, a short description of what they had read so far or which videos they had watched, with each programmatic assignment. Marks were ``perfect'' (for handing one in) or ``zero'' (for not handing one in). Of course this is a blatant bribe to get the students to read the textbook (or watch the course videos).  Many students initially thought of these as trivial ``free marks'' and of course they could use them in that way.  But the majority learned that these memos were a way to get detailed responses back, usually from the instructor or TA but sometimes from other students.  They learned that the more they put into a Reading Memo the more they got back.  The feedback to the instructor was also directly valuable for things like pacing.  Out of all the techniques we used, this one---the simplest---was the most valuable.

The other simple technique we used was discussion time.  Provocative, nearly paradoxical questions were the best for this.  For instance, consider the following classical paradox of the arrow, attributed to Zeno, interpreted in floating-point (actually, this was one of their exam questions this year):
\begin{lstlisting}{language=Matlab}
% Zeno's paradox, but updated for floating-point
%initial position of the arrow is s=0, target is s=1
s = 0
i = 0; % Number of times through the loop
% format hex
while s < 1,
    i = i+1;
    s = s + (1-s)/2 ;
end
fprintf( 'Arrow reached the target in %d steps\n', i)
\end{lstlisting}
In the original paradox, the arrow must first pass through the half-way point; and then the point half-way between there and the target; and so on, \emph{ad infinitum}.
The question for class discussion was, would the program terminate, and if so, what would it output?  Would it help to uncomment the \texttt{format hex} statement?  Students could (and did) type the program in and try it, in class; the results were quite surprising for the majority of the class.

\begin{figure}
\sidecaption
\centering
\includegraphics[width=7cm]{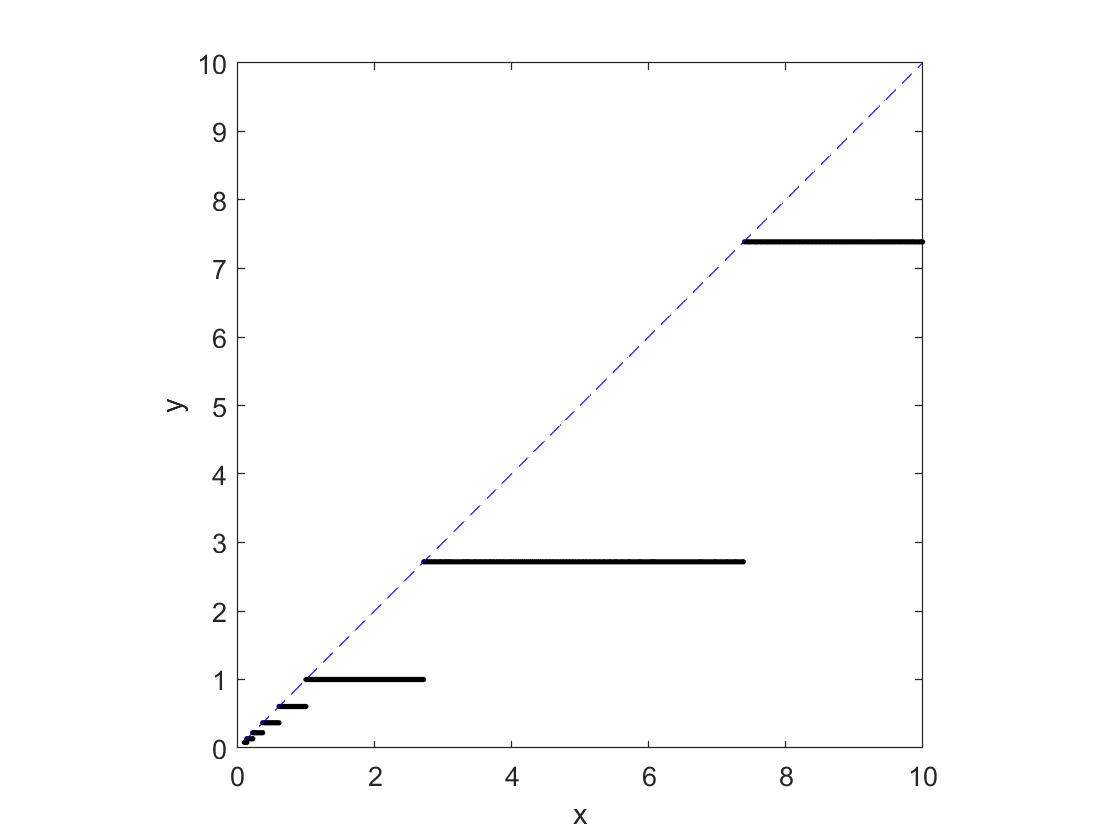}
\caption{The function $y = \textrm{Higham}(x) =  (x^{1/2^{52}})^{2^{52}}$, i.e. take $52$ square roots, and then square the result $52$ times, plotted for $2021$ points on $0.1 \le x \le 10$,
carried out in IEEE double precision.  Students are asked to identify the numerical values that $y$ takes on, and then to explain the result.  See Section 1.12.2 of~\cite{Higham(1996)} and also exercise 3.11 of that same book, and~\cite{kahan1980interval}. }\label{fig:Higham}
\end{figure}

Another lovely problem originates from one posed by Nick Higham: take an input number, $x$.  Take its square root, and then the square root of that, and so on $52$ times.  Now take the final result and square it.  Then square that, and again so on $52$ times.  One expects that we would simply return to $x$.  But (most of the time) we do \emph{not}, and instead return to another number.  By plotting the results for many $x$ on the interval $1/10 \le x \le 10$  (say) we see in figure~\ref{fig:Higham}, in fact, horizontal lines.  The students were asked to explain this.  This is not a trivial problem, and indeed in discussing this problem amongst the present authors, JHD was able to teach RMC (who has used this problem for years in class) something new about it.

We will \emph{not} give ``the answers'' to these questions here.  They are, after all, for discussion. [A useful hint for the repeated square root/repeated squaring one is to plot $\ln(y)$ against $x$.] We encourage you instead to try these examples in your favourite computer language and see what happens (so long as you are using floating-point arithmetic, or perhaps rounded rational arithmetic as in Derive)!

We will discuss, however, Kahan's proof of the impossibility of numerical integration~\cite{kahan1980handheld} here, as an instance of discussing the limits of technology.  This lesson must be done carefully: too much skepticism of numerical methods does much more harm than good, and before convincing students that they should be careful of the numerics they must believe that (at least sometimes) computation is very useful.  So, before we teach numerical integration we teach that symbolic integration is itself limited, especially if the vocabulary of functions is limited to elementary\footnote{The \emph{elementary functions} of the calculus are not ``elementary'' in the sense of being simple; but instead they are ``elementary'' in a similar sense to the elementary particles of physics.} antiderivatives.  As a simple instance, consider
\begin{equation}\label{eq:edelmanintegral}
  E = \int_1^\infty \frac{e^{-y^2/2}}{y+1}\>,
\end{equation}
which occurs in the study of the distribution of condition numbers of random matrices~\cite{Edelman1988}. The author laconically states that he ``knows of no simpler form'' for this integral.  In fact, neither do we, and neither do Maple or Mathematica: the indefinite integral is not only not elementary (provable by the methods of \cite{Davenport1986d}), it is right outside the reference books.  Of course the sad (?) fact is, as observed in~\cite{kahan1980handheld}, the vast majority of integrands that occur in ``real life'' must be dealt with numerically.  This motivates learning numerical quadrature methods.

However, it is a useful thing for a budding numerical analyst to learn that numerical techniques are not infallible, either. Consider the following very harmless function:  Aphra$(x) := 0$.  That is, whatever $x$ is input, the Aphra function returns $0$.  However, Aphra is named for \emph{Aphra Benn}, the celebrated playwright and spy for King Charles.  The function is written in Matlab in such a way as to \emph{record its inputs} $x$.

\begin{lstlisting}{language=Matlab, caption=A function named for Aphra Benn}
function [ y ] = Aphra( x )
%APHRA A harmless function, that simply returns 0
%
  global KingCharles;
  global KingCharlesIndex;
  n = length(x);
  KingCharles(KingCharlesIndex:KingCharlesIndex+n-1) = x(:);
  KingCharlesIndex = KingCharlesIndex + n;
  y = zeros(size(x));
end
\end{lstlisting}

If we ask Matlab's \emph{integral} command to find the area under the curve defined by Aphra$(x)$ on, say, $-1 \le x \le 1$, it very quickly returns the correct answer of zero.  However, now we introduce another function, called Benedict:

\begin{lstlisting}{language=Matlab}
function [ y ] = Benedict( x )
%BENEDICT Another harmless function
%   But this function is not zero.
  global KingCharles;
  global KingCharlesIndex;
  global Big;
  s = ones(size(x));
  for i=1:length(KingCharles),
      s = s.*(x-KingCharles(i)).^2;
  end
  y = Big*s;
end
\end{lstlisting}

\begin{figure}
\sidecaption
\centering
\includegraphics[width=7cm]{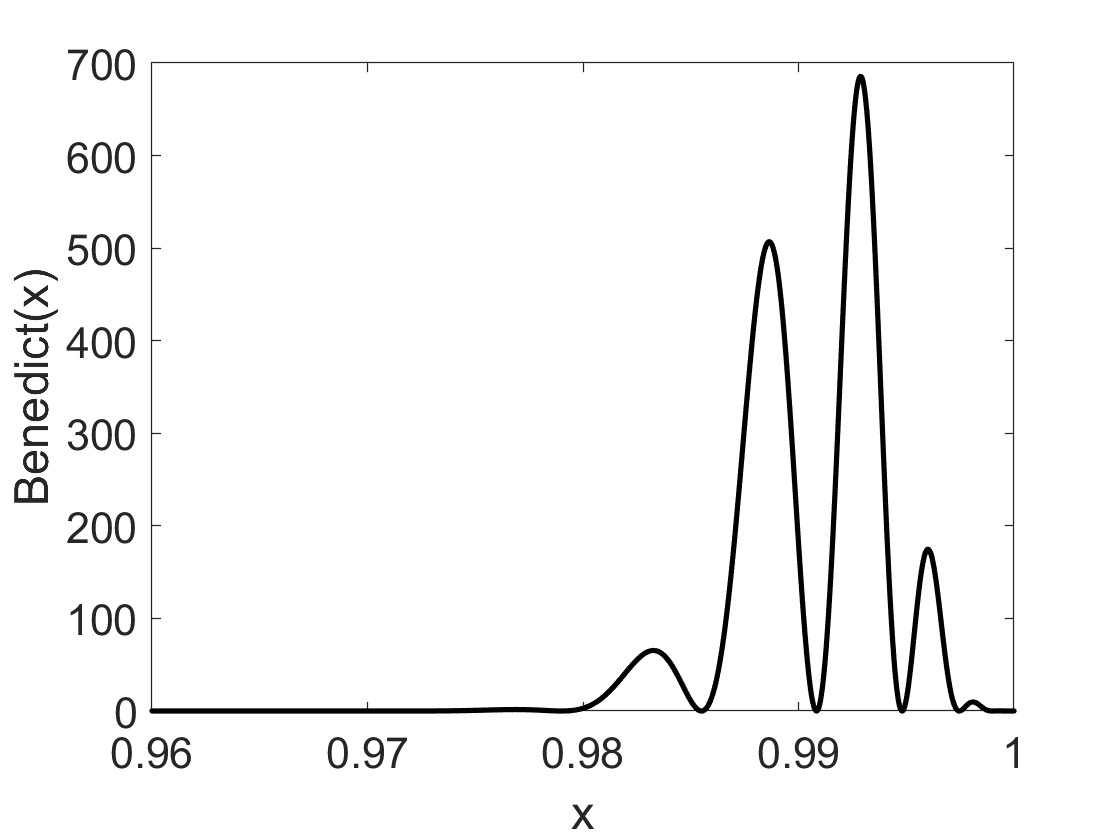}
\caption{The function Benedict$(x) = K\prod_{i=1}^{150} (x-s_i)^2$ where the $s_i$ are the $150$
sample points in the interval $-1 \le x \le 1$ reported by the function Aphra$(x)$, and with $K=10^{87}$.  We plot only an interesting region near the right endpoint of the interval.  We see that the area under this curve is not zero.}\label{fig:Benedict}
\end{figure}

This function is defined to be zero exactly at the points reported by Aphra, but strictly positive everywhere else: indeed the ``Big'' constant can be chosen arbitrarily large.  If we choose Big equal to $10^{87}$, then after calling Aphra with \verb+integral( @Aphra, -1, 1 )+ first, we find the function plotted in figure~\ref{fig:Benedict}.  It is clearly not zero, and indeed clearly has a positive area under the curve on the interval $-1 \le x \le 1$.

However, asking the Matlab built-in function \texttt{integral} to compute the area under Benedict$(x)$ on the interval $-1 \le x \le 1$ gives the incorrect result $0$ because the deterministic routine \texttt{integral} samples its functions adaptively but here by design the function Benedict traitorously behaves as if it was Aphra at the sample points (and only at the sample points).  This seems like cheating, but it really isn't: finding the good places to sample an integrand is remarkably difficult (and more so in higher dimensions).  One virtue of Kahan's impossibility proof is that it works for arbitrary deterministic numerical integration functions.  Without further assumptions (such as that the derivative of the integrand is bounded by a modest constant) numerical integration really is impossible.

The students \emph{do not like} this exercise.  They dislike learning that all that time they spent learning antidifferentiation tricks was largely wasted, and they dislike learning that computers can give wrong answers without warning.  Still, we feel that it is irresponsible to pretend otherwise.

Finally, the course was officially designated as an ``essay'' course.  This was in part recognition for the essay-like qualities of the lab reports, but was also explicitly in recognition of the similarities between a good computer program and a good essay: logical construction, clear division of labour, and good style.  It is our contention that programming and proving and explaining all share many attributes.  As Ambrose Bierce said, ``Good writing is clear thinking made visible.''

We also not only allowed but actively encouraged collaboration amongst the students.  They merely had to give credit to the other student group members who helped them, or to give us the name of the website they found their hints or answers on (frequently Stack Exchange but also Chegg and others).  Many students could not believe that they were being allowed to do this. The rationale is that in order to \emph{teach} something, the student had to know it very well.  By helping their fellow students, they were helping themselves more.

But modern programming or use of computers is \emph{not} individual heroic use: nearly everyone asks questions of the Web these days (indeed to answer some \LaTeX\ questions for the writing of this paper, we found the \LaTeX\ FaQ on the Help Page on Wikibooks useful; and this even though the authors of this present paper have \emph{decades} of \LaTeX\ experience).  We do not serve our students well if we blankly ban collaborative tools.  We feel that it is important to teach our students to properly \emph{acknowledge} aid, as part of modern scientific practice.

\subsection{Assessment}
But we did not allow collaboration on the midterm exam, which tested the students' individual use of Matlab at computers locked-down so that only Matlab (and its help system) could be used.  Examination is already stressful: an exam where the student is at the mercy of computer failure or of trivial syntax errors is quite a bit more stressful yet.  To mitigate this we gave \emph{practice exams} (a disguised form of active learning) which were quite similar to the actual exam.  The students were grateful for the practice exams, \emph{and moreover found them to be useful methods to learn}.

Exam stress---assessment stress in general---unfortunately seems to be necessary\footnote{Given the economic constraints of the large class model, we mean. Even then, there may be alternatives, such as so-called ``mastery grading''~\cite{armacost2003using}. We look forward to trying that out. Exam stress is often counterproductive, and the current university assessment structures do encourage and reward successful cheating.  We would like a way out of this, especially now in COVID times.}: if the students \emph{could} pass the course without learning to program Matlab, they \emph{would} do so, and thereafter hope that for the rest of their lives they could get other people to do the programming.  Students are being rational, here: if they were only assessed on mathematical knowledge and not on programming, then they should study the mathematics and leave the programming for another day.  So we must assess their individual programming prowess.

In contrast, the students were greatly relieved to have a final exam that ``merely'' asked them to (in part) write pencil-and-paper programs for the instructor to read and grade.  In that case, trivial errors---which could derail a machine exam---could be excused.  On the other hand, the instructor could (and did) ask for explanations of results, not merely for recitations of ways to produce them.
%\section{Introduction to Computer Algebra}
%\label{sec:computeralgebra}

\section{Computational Discovery/Experimental Mathematics}
\label{sec:computationaldiscovery}
The courses that we describe in this section are described more fully elsewhere~\cite{chan2017random,Chan:2022:CDJ}.  Here we only sketch the outlines and talk about the use of active learning techniques with (generally) introverted mathematics students.  The major purpose of these related courses (a first-year course and a graduate course, both in Experimental Mathematics, taught together) was to bring the students as quickly as possible to the forefront of mathematics.

\begin{quotation}
  Short is the distance between the elementary and the most sophisticated results, which brings rank beginners close to certain current concerns of the specialists.\\
\hfill---\cite{mandelbrot2002some}
\end{quotation}

In this we were successful.  For example, one student solved a problem that was believed at the time to be open (and she actually solved it \emph{in-class}); although we were unaware at the time, it turned out to have actually been solved previously and published in 2012, but nonetheless we were able to get a further publication out of it, namely~\cite{li2019revisiting}, having taken the solution further.  There were other successes.  Some of the projects became Masters' theses, and led to further publications such as~\cite{chan2017new}, for example.

The course was also \emph{visually} successful: the students generated many publication quality images,
some of which were from new \href{bohemianmatrices.com}{Bohemian Matrix} classes.  Indeed some of the images at that website were produced by students in the course.

\subsection{Choice of Programming Language}
\begin{figure}
\sidecaption
\centering
\includegraphics[width=7cm]{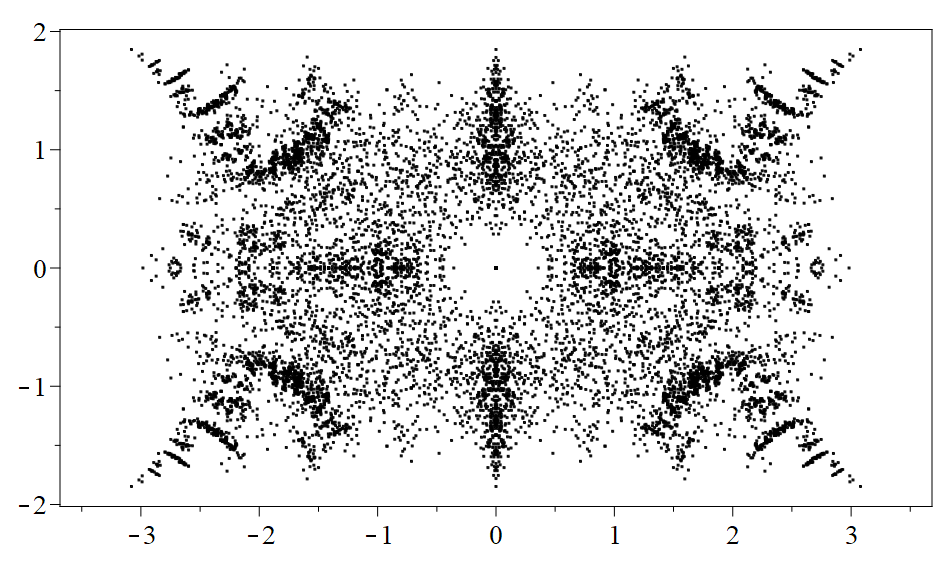}
\caption{All the complex eigenvalues of all the $7$-dimensional symmetric tridiagonal (but with zero diagonal) matrices with population $\{-5/3-i, -5/3+i, 5/3+i, 5/3-i\}$.  There are $4^6=4096$ such  matrices, but only about half as many distinct characteristic polynomials in the set. }\label{fig:ElectricBohemian}
\end{figure}

We used Maple for this course, because its symbolic, numerical, and visual tools make it eminently suited to experimental mathematics and computational discovery; because it was free for the students (Western has a site licence); and because of instructor expertise~\cite{corless2004essential}.  For instance, Maple allows us to produce the plot shown in figure~\ref{fig:ElectricBohemian} of all the eigenvalues of a particular class of matrices.  This figure resembles others produced by students in the class, but we made this one specifically for this paper.  There are $4096$ matrices in this set, each of dimension $7$.  However, there are only $2038$ distinct characteristic polynomials of these matrices because some are repeated.  Getting the students to try to answer questions such as ``how many distinct eigenvalues are there'' is a beginning (this is not obvious, because again there are repeats: the only way we know how to answer this is to compute the GCD-free basis of the set of $2038$ degree $7$ polynomials, in fact).  A bigger goal---in fact, the main goal of the course---was getting the students to come up with their own questions.  It helped that the students were encouraged to invent their own classes of matrices (and they came up with some quite remarkably imaginative ones).

\subsection{Pedagogical Methods}
This course was designed wholly with active learning in mind.  It took place in the Western Active Learning Space, which was divided into six tables called Pods, each of which could seat about seven students; the tables were equipped with technology which allowed students to wirelessly use the common screens to display materials to each other.  The smartboards were (in principle) usable in quite sophisticated technological ways; in practice, the varieties of whiteboards with simple coloured pens were just as useful.

Students enrolled in the first-year course were grouped with students enrolled in the graduate course.  Each group benefitted from the presence of the other: the presence of the senior students was a calming factor, while the junior students provided significant amounts of energy.  The grad student course also had an extra lecture hour per week where more advanced topics were covered in a lecture format.

Active learning techniques run from the obvious (get students to choose their own examples, and share) through the eccentric (interrupt students while programming similar but different programs and have them trade computers and problems) to the flaky (get them to do an interpretive dance or improvisational skit about their question). We tried to avoid the extremely flaky, but we did mention such, so that these introverted science students knew that this was within the realm of possibility.

The simplest activity was typing Maple programs that were handwritten on a whiteboard into a computer: this was simple but helpful because students learned the importance of precision, and had \textsl{immediate} help from their fellow students and from the TA.

Next in complexity was interactive programming exercises (integrated into the problems). Mathematicians tend to under-value the difficulty of learning syntax and semantics simultaneously.  The amplification of human intelligence by coupling it with computer algebra tools was a central aspect of this course.

We describe our one foray into eccentricity. The paper Strange Series and High Precision Fraud by Borwein and Borwein \cite{borwein1992strange} has six similar sums. We had six teams program each sum, at a stage in their learning where this was difficult (closer to the start of the course). After letting the teams work for twenty minutes, we forced one member of each team to join a new team; each team had to explain their program (none were working at this stage) to the new member. This exercise was most instructive. The lessons learned included:
\begin{itemize}
	\item people approach similar problems very differently
	\item explaining what you are doing is as hard as doing it (maybe harder)
	\item basic software engineering (good variable names, clear structure, economy of thought) is important
	\item designing on paper first might be a good idea (nobody believed this, really, even after)
	\item social skills matter (including listening skills).
\end{itemize}
\subsection{Assessment}
The students were assessed in part \emph{by each other}: we used peer assessment on class presentations.  The instructor informed the students that he would take their assessments and \emph{average them with his own} because peer assessment is frequently too harsh on other students; they found this reassuring.  The main mark was on an individual project, which took the full term to complete.  They had to present intermediate progress at a little past the half-way point.  Marks were also given for class participation.

Collaboration was encouraged.  The students merely had to make proper academic attribution.  While, technically, cheating might have been possible---one might imagine a plagiarized project---there was absolutely no difficulty in practice.  The students were extremely pleased to be treated as honourable academics.

\section{Programming and Discrete Mathematics}
This course described in this section is also more fully explained elsewhere: see~\cite{betteridge2019teaching}. We restrict ourselves here to an outline of the goals and methods.

The course XX10190, Programming and Discrete Mathematics, at the University of Bath is both similar and dissimilar to the Western courses described above.  One of the big differences is that it was designed specifically for the purpose of teaching programming to mathematical scientists by using mathematics as the proving ground.  The course was designed after significant consultation and a Whole Course Review in 2008/2009.  In contrast, the Western course designs were driven mostly by the individual vision of the instructor.  The Bath course therefore has a larger base of support and is moreover supported by the recommendation from~\cite{bond2018era} that ``every mathematician learn to program.''  As such, it is much more likely to have a long lifespan and to influence more than a few cohorts of students; indeed, since it has been running for ten years, it already has\footnote{One British citizen in 25,000 is a graduate of XX10190.}. Now that RMC has retired from Western and the numerical analysis course has been taken over by a different instructor, the course there is already different. Conversely all the Bath authors have moved on from XX10190, but the course is much the same. This is the differential effect of institutional memory.

Another big difference is that the course is in first year, not second year; more, it runs throughout the first year, instead of only being a thirteen week course.  This gives significant scope for its integrated curriculum, and significant time for the students to absorb the lessons.

However, there are similarities.  The focus on discrete mathematics makes it similar to the Experimental Mathematics courses discussed above, with respect to the flavour of mathematics.  Indeed, perhaps the text~\cite{eilers2017introduction} might contain some topics of interest for the course at Bath.  Although the focus is on discrete mathematics, some floating-point topics are covered and so the course is similar to the Numerical Analysis course above as well.  But the main similarity is the overall goal: to use mathematical topics to teach programming to mathematical scientists, and simultaneously to use programming to teach mathematics to the same students.  This synergistic goal is eminently practical: learning to program is an effective way to learn to do mathematics.

Another similarity is respect for the practical \emph{craft} of programming: the papers~\cite{davenport2016innovative} and~\cite{wilson2006software} discuss this further. To this end, the instructors use Live Programming~\cite{rubin2013effectiveness}, defined in~\cite{paxton2002live} as ``the process of designing and implementing
a [coding] project in front of class during lecture
period.'' This is in contrast to the Western courses, where an accidental discovery proved valuable: the instructor was for several years discouraged from using keyboards owing to a repetitive strain injury, and as a consequence took to writing code on the whiteboard.  This had unexpected benefits when the students would ask him to debug their code, and he would do so in a Socratic manner by asking the students to relay error messages.  In doing so, the students frequently found their own solutions.  However, one of the most common requests from students was for live demonstrations: there is no question that live programming techniques can be valuable.

\subsection{Choice of Programming Language}
A major similarity to the Western course is the choice of programming language: Matlab.  As with the Western course, Matlab may eventually be displaced by Python, but is an admirable first language to learn for mathematical scientists.  This choice came with several unanticipated benefits, as described in~\cite{betteridge2019teaching}: for instance, the statisticians teaching subsequent courses found it simpler to teach R to students who had a working knowledge of the similar language Matlab.

\subsection{Pedagogical Methods}
The course is fifty percent Programming and fifty percent Discrete Mathematics.  The course is team taught, with Lecturers and Tutors.  The whole cohort have one programming lecture, one Discrete Mathematics lecture, and one Examples class per week. The roughly $300$ students are divided up into tutorial groups of size $20$, and there is one Discrete Math tutorial per week (when budgets allow: some years this has been financially impossible) and one Programming Lab on Fridays, after the whole-cohort classes (this apparently minor timetabling point is pedagogically very helpful).  Management of this relatively large staff with its hierarchical structure repays attention, and the instructors have found it valuable to provide tools such as a separate mailing list for tutors.  The course uses Moodle and its various electronic delivery tools.

The Lab physically holds 75 seats, divided into five tables with fifteen computers each.  There is one tutor for approximately ten students: students and tutors are assigned to specific groups of seats. This division allows greater and more sustained personal contact, and more active learning.

Tutors must take great care helping students in labs.
The student is not just learning a language but a new logical structure, while instructors are proficient coders.
When a student asks for help, it is far too easy for a tutor to `fix' the code for them, particularly when one is new to teaching.
While this is the path of least resistance, because the student's priority is working code, for many not only does this do little for learning but in fact this can be detrimental to learning.
If a tutor rewrites code with no sympathy for the student's approach, this can just alienate and destroy confidence.

A philosophy of `never touch the keyboard' embodies our approach. As one practices, this approach reveals subtler layers.
[We have also noted that with remote teaching, although one is physically removed, practising the method is more difficult!]
The philosophy applies to both instructor and student. It really means not telling students the difficulty with their draft code, but rather discovering it with them.
One method is to ask what the student is trying to do, read their code with them, and try to
nurture \emph{their} creativity.
It can be time intensive, and is not easy. One needs react to the student, taking care not to add to the student's pain by repeating the same question\footnote{Although it's true that, sometimes, simply reading a question aloud can be surprisingly useful; but of course tone matters, here. Reading the question aloud as if it were a reminder to the \emph{instructor} can be less painful for the student.}, methods like pseudocode and flow diagrams can be useful for withdrawing from the screen.
Any suffering (on both sides) is justified when the students `gets it' and the sparks of understanding light in their eyes.

\subsection{Assessment}
Similar to the ``Reading Memos'' of the Western courses, the Bath course has what is called a ``tickable.''  These are short exercises---gradually increasing in difficulty throughout the year---which are graded only on a Yes/No basis. A tickable therefore differs from a Reading Memo in that it requests some well-defined activity, whereas a Reading Memo is less well-defined and more open-ended.  The similarity is in their assessment and in their encouragement of continual work throughout the course.

For instance one tickable from midway through the first semester is given here:
\medskip\par\noindent
\textbf{Tickable:} Write a recursive Matlab function, in the file \verb`myexpt.m`, which will compute $A^n$ (via the call \verb`myexpt(A,n)`) using equation (\ref{eq:fastexp}), for any square matrix $A$.
\begin{equation}
	x^n =
    \begin{cases}
    	1 & \text{ if } n = 0\\
        (x \cdot x)^{n/2} & \text{ if } n \text{ is even}\\
        x \cdot (x \cdot x)^{(n-1)/2} & \text{ if } n \text{ is odd}
    \end{cases}
    \label{eq:fastexp}
\end{equation}

This tickable is then used to write another small program for quickly calculating the $n^\text{th}$ Fibonacci number. During lab sessions, a tutor (who has approximately 7--10 students assigned to be their tutees for the whole semester, or ideally year) walks around the computer terminals offering help with the mathematical or programming aspects of the exercise. Students who successfully get this code running can also re-use this routine for parts of the coursework at the end of the semeseter.

An insufficient number (fewer than $80$\% of the total) of tickables marked ``Yes'' results in a pro rata  reduction in the summative mark.  This is widely perceived as fair, because there is general agreement that doing the work as you go along helps you to learn the material.

Otherwise there is significant use of automatic assessment tools via Moodle, with tutors providing more detailed feedback on programming style.

\section{Outcomes}
In both the Western and the Bath cases, the student surveys showed great satisfaction.  For instance, the TA for the Western Numerical Analysis course twice won the ``TA of the Year'' award from the Society of grad students.  True measurement of the effectiveness of these courses is naturally difficult, but the indications pointed out in~\cite{betteridge2019teaching}, which include superior outcomes in downstream courses, seem quite solid.

Since no controlled experiments were made about teaching methods---in neither case was there a control group, where different methods were used---this kind of qualitative good feeling about outcomes may be the best indication of success that we can obtain.  This clearly touches on the ethics of testing different methods of teaching, and we take this up briefly in the next section.

\section{Ethics, Teaching, and Eudaemonia}
Much published research on teaching just barely skirts rules about experimentation on humans. The `out' that is most frequently used is the \emph{belief} on the part of the teachers that what they are doing is ``best practice''.  It is rare to have proper statistical design with control groups to compare the effects of an innovation with mere placebo change over the status quo.  The previously mentioned research on Active Learning includes some that meets this stringent standard, and as previously mentioned the evidence is so strong that it is now known to be \emph{unethical} not to use active learning.  Still, active learning is labour-intensive (on everyone's part---it's a lot simpler for a student to sit and pretend to listen in class, and then cram for an exam in the traditional ``academic bulimia'' model) and not everyone is willing to pay the price for good ethics.

Another significant piece of active learning is the social aspect.  Humans are social animals and teaching and learning is part of how we interact in person.  University students appear to value \emph{personal contact} above nearly anything else~\cite{seymour1997talking}.  Working against that, economics of scale mean that universities want to provide certificates of learning by using only small numbers of teachers for many students; this impersonal model is already unsatisfactory for many students.  This time of social isolation due to COVID-19 is making this worse, of course, in part because teaching and learning are becoming even more impersonal.  One response to this pressure---and this was happening before COVID---is to try to let computers help, and to use AI to personalize instruction and especially assessment.

There is an even deeper ethical question at work, however.  A teacher who taught lies\footnote{Except as an important stepping stone to the real truth---see the entry ``Lies to Children'' in Wikipedia.  Sometimes a simplistic story is the right first step.} would be properly viewed as being unethical, even as being evil. A teacher who hid important facts from the students would be scarcely less unethical. This observation seems to be culturally universal (with perhaps some exceptions, where secret knowledge was jealously guarded, but valued all the more because its exclusiveness).  Yet, aside from idealism, what are the motivations for the teacher to tell the truth, the whole truth, and nothing but the truth?

When humans are the teachers, this is one question.  We collectively know to be skeptical of the motives of people: who benefits from this action, and why are they doing this?  Teaching pays, and not only in money: perhaps the most important part of our pay is the respect of those that we respect.  Most of us understand that the best teachers do their jobs for love of watching their students understand, especially seeing ``light-bulb moments''.  But when the teacher is an app on your smartphone, the questions become different.  We will take as example the popular language app Duolingo~\cite{von2013duolingo}. The goals of a company that sells (or gives away---Duolingo is free by default, supported by advertising) an app to teach you something may very well be different to the goals of a human teacher.  Indeed, and there is nothing hidden or nefarious about this, one of the goals of the maker of Duolingo is to \emph{provide low-cost translation services}, essentially by distributing the translation tasks to (relatively) trusted human computers.  It is an ingenious idea: make the skilled app user pay for the service of learning a new language by providing some services, more as the learning progresses, that others want.  The question then becomes not ``what does my teacher gain'' but rather ``what does the creator of this service gain;'' more insidiously, if a teaching app became truly viral, it might be ``what reproductive value does this app gain?''

The modern university system has evolved from its religious roots to provide a desired service today---namely access to the scholarship of the world---to anyone who can find a way to access the University system.  We (mostly) share a belief that access to education is one of the great benefits, and provides the key to a better life, a good life, the best life possible (indeed to \emph{eudaemonia}, in Aristotle's term\footnote{Aristotle may have done us a disservice by looking down on crafts and craftspeople; the term Software Carpentry is not likely to induce respect for the discipline in academia, for instance.  We lament this prejudice.}, although people still argue about what exactly he meant by that). It is not at all clear to us that an artificially intelligent teacher (even if humans are in the loop, as with Duolingo) would necessarily share this belief.  The benefits to such a ``teacher'' of actively \emph{discouraging} critical thinking are unfortunately extremely clear: one only has to look at the unbearable success of lies on social media to see the problem.

It seems clear to us that we as teachers should pay attention to the ethics of teaching by or with the help of AIs.

\section{Concluding Remarks}

\begin{quotation}
  Instead there must be a more serious concern with the significant ways in which computational resources can be used to improve not so much the \textbf{delivery} but rather the \textbf{content} of university courses.
  \begin{flushright}---\cite{abelson1976computation}\end{flushright}
\end{quotation}

%\begin{quotation}
%Please do not use quotation marks when quoting texts! Simply use the \verb|quotation| environment -- it will automatically be rendered in line with the preferred layout.
%\end{quotation}
The content of mathematics courses has changed over the past few decades (this has been noted in many places, but see e.g.~\cite{corless1997scientific}).  Some of that change has been forced by the increasing number of students and their necessarily more diverse backgrounds and interests; some of that change has been deliberate abandonment of no-longer-useful techniques; and some of that change has been driven by the introduction of new tools.  One new tool that we have not yet talked about is Wolfram Alpha.  This is nearly universally available, free, almost uses natural language input---it's pretty accepting, and the students find it simple to use---and produces for the most part very legible, neat, and useful answers to problems at roughly the first year university level.  We believe that its use (or the use of similar tools) should not only be allowed in class, but encouraged.  The students will still be \emph{responsible} for the answers, and it helps to give examples where Wolfram Alpha's answers are wrong or not very useful; but it is irresponsible of us to ignore it.  Matlab, Maple, Mathematica, Python, NumPy, and SymPy provide other tools for mathematical thinking, on a larger scale.  We believe that it is incumbent on us as educators to teach our students the kinds of mathematics that they can do when using those tools.

The requirement for correctness and reproducibility in mathematical sciences is paramount, but academia has been slow to apply this rigorously to its codes.
In software development this is achieved with testing, validation and version control.
While comparison with expectation and (better) analytic proof are adequate for validation we have not formally taught testing or version control in our undergraduate programmes.
The time pressure on curriculum cannot excuse this omission much longer.
The value of adopting professional practices goes beyond those who will work as software engineers.
They are vital tools for working efficiently, contributing to open software,  for data scientists and AI engineers to manage data and to ensure trust in the methods that they develop and apply in their careers.
These enable students to use their computational tools responsibly.

We have not talked here about GeoGebra, which is probably now the most popular computational thinking tool for mathematics in the world.  This is because we are ``old guard'' (well, some of us are) and GeoGebra is a newer tool, one that we have not yet used.  However, it seems clear to us that the same principles that we have been using for our other tools also apply here: the students should be aware of the program's limitations; the students should know when the answer is correct and when it is not; and the students should be responsible for the answers.

\begin{flushright}
Plus \c{c}a change, plus c'est la m\^eme chose. ---Jean Baptiste Alphonse Karr, 1849
\end{flushright}

With the advent of modern AI tools for education, more questions arise.  We believe that amongst the most important questions for AIEd will be about the \emph{ethics} underlying the tools.  We all know now that machine learning can easily copy our biases and prejudices, without us intending; we also know that the goals of developers of AIEd tools may well be different than the goals of good teachers\footnote{See also \cite{Bradfordetal2009a}, which shows that this can affect the basic meaning of equality: pedagogical equality is not the same as mathematical equality.  It is perfectly possible for two expressions to be mathematically equal, but only one expression to be the desired student response.}.  The ethics of AIEd is beginning to be studied intensively (see e.g.~\cite{aiken2000ethical,sijing2018artificial}) but clearly we are only just scratching the surface of the issues, which include some very deep classical philosophical problems, including how to live a good life (achieve eudaemonia).  The amplified human life, when humans use computers to increase their thinking capability, clearly also needs philosophical study.  Not only philosophers, but cognitive scientists, as well as computer scientist experts in AI, will be needed to properly develop these tools.

\begin{acknowledgement}
RMC thanks the Isaac Newton Institute for Mathematical
Sciences and the staff of both the University Library and the Betty and Gordon Moore Library at Cambridge for support and hospitality during the programme
Complex Analysis: Tools, techniques, and applications, by EPSRC Grant \# EP/R014604/1
when some of the work on this project was undertaken. RMC likewise thanks the University of Bath for an invitation to visit Bath, at which this project was started.
EYSC and RMC also thank Western University for a grant to work on the project \emph{Computational Discovery on Jupyter}, some of whose results are discussed here.
\end{acknowledgement}
%
%\section*{Appendix}
%\addcontentsline{toc}{section}{Appendix}
%
%
%When placed at the end of a chapter or contribution (as opposed to at the end of the book), the numbering of tables, figures, and equations in the appendix section continues on from that in the main text. Hence please \textit{do not} use the \verb|appendix| command when writing an appendix at the end of your chapter or contribution. If there is only one the appendix is designated ``Appendix'', or ``Appendix 1'', or ``Appendix 2'', etc. if there is more than one.

\bibliographystyle{apalike}
\bibliography{teaching}
\end{document}